\documentclass[11pt]{amsart}
\usepackage[utf8]{inputenc}
\usepackage{amsmath, palatino, mathpazo, amsfonts, amssymb, mathtools}
\usepackage[mathscr]{eucal}
\usepackage[all]{xy}
\usepackage{datetime}
\usepackage{amsthm}
\usepackage{comment}
\usepackage{marginnote}
\usepackage{tikz}
\usepackage{pgfplots}
\usepackage{subcaption}

\theoremstyle{definition}
\newtheorem{theorem}{Theorem}[section]

\newtheorem{proposition}[theorem]{Proposition}

\newtheorem{definition}[theorem]{Definition}

\theoremstyle{remark}
\newtheorem{conjecture}[theorem]{Conjecture}

\numberwithin{equation}{section}
\allowdisplaybreaks

\title{A simple proof of Dahmen's conjectures}

\author{Po-Sheng Wu}
\date{\today}

\begin{document}

\maketitle
\begin{abstract}
The number of Lame equations with finite (ordinary or projective) monodromy has been conjectured by S. R. Dahmen, and a few proofs have been proposed. It is known that Lame equations with unitary monodromy are corresponding to spherical tori with one conical singularity, and the geometry of such surfaces had been studied with triangulation recently. In this paper, we will apply the results on spherical tori to give an alternative proof of Dahmen's conjectures.
\end{abstract}

\section{Introduction}

Given a lattice $\Lambda=\mathbb Z\omega_1+\mathbb Z\omega_2$ on $\mathbb C$ with $(\omega_1,\omega_2)=(1,\tau), \textnormal{Im}(\tau)>0$, the Lam\'e equation on the elliptic curve $E=\mathbb C/\Lambda$ is a second order ordinary differential equation:
\begin{equation}\label{lame}
 \frac{\partial^2 w}{\partial z^2} - \Big( n(n+1) \wp(z)+B \Big) w =0,
\end{equation}
where $B\in \mathbb C$ and $\wp$ is the Weierstrass elliptic function. We consider the case $n\in \mathbb Z_{>0}$ and study number of Lam\'e equations with given finite monodromy groups. \\ 

It is known that all the finite monodromy groups of Lam\'e equations are cyclic. The following two conjectures are proposed by S.R. Dahmen, and later he proved the first conjecture using dessin d'enfant \cite{S1}\cite{S2}.

\conjecture Let $L_n(N)$ be the number of Lam\'e equations (\ref{lame}) with projective monodromy group isomorphic to the cyclic group $C_N$, then
\begin{equation}\label{fin_proj_mono}
L_n(N)=\dfrac{n(n+1)}{12}(\Psi(N)-3\phi(N))+\dfrac{2}{3}\epsilon(n,N),
\end{equation}
for $N\geq 3$, where \[\epsilon(n,N)=\begin{cases}1,\quad \textnormal{if }n=3\textnormal{ and }3|N-1,\\0,\quad \textnormal{otherwise.}\end{cases},\] $\phi$ is the Euler totient function and \[\Psi(N)=|\{(k_1,k_2)\in \{0,\dots,N-1\}^2, \textrm{gcd}(k_1,k_2,N)=1\}|\].\\

\conjecture Let $L'_n(N)$ be the number of Lam\'e equations (\ref{lame}) with ordinary monodromy group isomorphic to the cyclic group $C_N$, then
\begin{equation}\label{fin_ordi_mono}
L'_n(N)=\dfrac{1}{2}(\dfrac{n(n+1)}{24}\Psi(N)-a_n\phi(N)-b_n\phi(N/2))+\dfrac{2}{3}\epsilon(n,N),
\end{equation}
for $N\geq 3$, where $a_{2l}=a_{2l+1}=l(l+1)/2$, $b_{2l-1}=b_{2l}=l^2$, and $\phi(N/2)=0$ if $N$ is odd. \\

Note that in Dahmen's paper, the conjectures were originally stated for pushforward Lam\'e equation on $\mathbb{CP}^1$, thus the finite monodromy groups they concerned becomes the dihedral groups $D_N$ instead of $C_N$.\\

Up to $N\leq 4$ the second conjecture has been verified ($n=1,2,3$ by S.R. Dahmen, and $n=4$ by Y.C. Chou via modular form calculations (\cite{MF2}, appendix)). On the other hand, Z. Chen, T.J. Kuo and C.S. Lin have announced a proof of the second conjecture \cite{PE}, which involves more technical analyses on Painlev\'e equations. In this paper we give a simple and uniform proof of these two conjectures from the perspective of spherical tori, developed by A. Eremenko et al.\ \cite{E}.\\

In Section 2 we will briefly introduce results known from literatures, and we will complete our proof of main theorems in Section 3.

\section*{Acknowledgement}
This paper is adapted from part of my 2021 master thesis \cite{TH}. I would like to appreciate Prof. Chin-Lung Wang for his guidance on this vast topic. Thanks to Prof. Hui-Wen Lin for encouraging me to join this project. Thanks also to You-Cheng Chou, Hung-Hsun Yu and Ting-Wei Chao for discussing and reviewing my paper. I am also grateful to my family for their financial and mental support.

\section{Known results}
We first recall some basic results on Lam\'e equations with unitary monodromies. The details can be found in \cite{WW}\cite{A}. It is a classical result that the Lam\'e equation (\ref{lame}) has the following ansatz solution
\proposition \label{ansatz}
\begin{equation}
w_{\underline a}=\exp(z\sum_{\mu=1}^n\zeta(a_\mu))\prod_{\mu=1}^n\frac{\sigma(z-a_\mu)}{\sigma(z)}
\end{equation}
where $\underline a=(a_1,\dots,a_n)$ satisfies \[\sum_{\nu\neq \mu} (\zeta(a_\nu)-\zeta(a_\mu)+\zeta(a_\mu-a_\nu))=0\] for $\mu=1,\dots,n$, and $B=(2n-1)\sum_{\mu=1}^n\wp(a_i)$.\\ 

Another solution for the same equation (with same $B$) can be chosen as $w_{-\underline a}$, where $-\underline a=(-a_1,\dots,-a_n)$, except that for $2n+1$ values of $B$ such that $\underline a=-\underline a$ up to a permutation. For ansatz solutions with $\underline a\neq-\underline a$, the monodromy of the quotient $f=\omega_{\underline a}/\omega_{-\underline a}$ (i.e. the projective monodromy of the equation) is given by
\proposition
\begin{equation}
f(z+\omega_i)=\exp(\int_{\gamma_i} g)f(z),\qquad i=1,2,
\end{equation}
where \[g=(\log f)'=\displaystyle\sum_{\mu=1}^n \dfrac{\wp'(a_\mu)}{\wp(z)-\wp(a_\mu)}\] and $\gamma_1,\gamma_2$ are the two fundamental loops on the torus.\\

Also note that the ansatz gives trivial monodromy at the singularity at $0$. Thus, the condition of a Lam\'e equation having unitary monodromy group is equivalent to that \[\int_{\gamma_i} g\in \sqrt{-1}\mathbb R,  i=1,2.\] 

A little reduction using Legendre's relation shows that 

\proposition (\cite{A}, p5-8, p21) The unitary monodromy condition for Lam\'e equations is equivalent to that \begin{equation}\label{green-hecke}
\sum_{\mu=1}^n Z(a_\mu)=0\end{equation}

Here if we write $a_\mu=t_\mu\omega_1+s_\mu\omega_2$ with $s_\mu,t_\mu\in \mathbb R$, and the quasi-periods $\eta_j=\zeta(z+\omega_j)-\zeta(z)$, then $Z$ is the Hecke function defined as \[Z(a_\mu)=\zeta(a_\mu)-t_\mu\eta_1-s_\mu\eta_2\]

If we further denote $s=\sum_{\mu=1}^n s_\mu$ and $t=\sum_{\mu=1}^n t_\mu$, then the projective monodromy is given by  \begin{equation}\label{proj_mono_eq}
\begin{cases}
f(z+\omega_1)=\exp(-4i\pi s)f(z),\\ 
f(z+\omega_2)=\exp(4i\pi t)f(z).
\end{cases}
\end{equation}
while the ansatz $w_{\underline a},w_{-\underline a}$ have monodromy
\begin{equation}\label{nonproj_mono_eq}
\begin{cases}
w_{\pm \underline a}(z+\omega_1)=\exp(\mp 2i\pi s)w_{\pm \underline a}(z),\\ w_{\pm \underline a}(z+\omega_2)=\exp(\pm 2i\pi t)w_{\pm \underline a}(z).
\end{cases}
\end{equation}

Thus $s,t (\mathrm{mod }1)$ determines the ordinary monodromy of the equation, while $2s,2t (\mathrm{mod }1)$ determines the projective monodromy.\\

The quotient $f$ can be viewed as the developing map of a spherical torus with one conical singularity.

\definition A spherical torus (with one conical singularity, omitted for short) $(S,x)$  is an oriented Riemannian surface $S$ of constant curvature $1$ and genus $1$, with a conical singularity $x$ of angle $2\pi\theta$, i.e., there is an local isometry of $S$ at $x$ to a spherical fan with corner of angle $2\pi\theta$ identifying its two edges.\\

We can also give any spherical torus $S$ a complex structure by identifying $S^2\cong\mathbb {CP}^1$, and this will make $S$ a Riemann surface of genus $1$ with a puncture at $x$.\\

\proposition If equation (\ref{lame}) has unitary monodromy for $f=w_1/w_2$, then the pullback of the Fubini-Study metric on $\mathbb{CP}^1\cong S^2$ by $f$ will produce a spherical torus $(E,x)$ with a conical singularity of angle $(4n+2)\pi$, and $E$ is its underlying Riemann surface. Conversely, any spherical torus of angle $(4n+2)\pi$ arises from a Lam\'e equation (\ref{lame}) with unitary monodromy, with the underlying elliptic curve given by the identification above.\\

Note that $f$ is depending on the choice of $w_1$ and $w_2$, so there is a correspondence between Lam\'e equation with unitary monodromy, and projective equivalence classes of spherical surfaces.

\definition We say two spherical tori are projective equivalent if their developing maps $f_1,f_2$ are differ by a composition of m\"obius transformation $\gamma\in\textrm{PSL}(2,\mathbb C)$ on $S^2\cong \mathbb{CP}^1$, i.e., $f_2=\gamma\circ f_1$, or equivalently, $f_1$ and $f_2$ correspond to the same Lam\'e equation.\\

In the following we summarize the results on spherical tori from \cite{E} which are needed in the paper.

\definition A spherical triangle is an oriented Riemannian surface $P$ of constant curvature $1$ with three geodesic boundaries.\\

Although the case $\theta\not\in 2\mathbb Z+1$ is not used in our paper, we include the result for completeness.

\proposition (Theorem B in \cite{E}) Let $(S,x)$ be a spherical torus. If $\theta\not\in 2\mathbb Z+1$, then $S$ can be decomposed into two isometric spherical triangles, with the three interior angles $\pi\theta_1,\pi\theta_2,\pi\theta_3$ satisfying the triangle inequalities, i.e., $|\theta_1-\theta_2|\leq\theta_3\leq\theta_1+\theta_2$. Conversely, any such spherical triangle will uniquely determine a spherical torus, except for the case when some of $\theta_i=\theta_j+\theta_k$ holds, $\{i,j,k\}=\{1,2,3\}$, in which case the spherical triangle and its mirror image correspond to the same spherical torus.\\

For any spherical torus with such decomposition, we denote by $\triangle,\triangle'$ the two spherical triangles, and $L_i,P_i$ ($L_i',P_i'$) the edges and vertices on $\triangle$ ($\triangle'$, respectively), so that $P_1,P_2,P_3$ are ordered clockwise on $\partial\triangle$. To glue $\triangle$ and $\triangle'$ into $(S,x)$, we simply glue $L_i$ along $L_i'$ so that the orientations on both sides are compatible, $i=1,2,3$. The vertices of the spherical triangles would then form a conical singularity.

\proposition (Theorem E in \cite{E}, with a $2$-torsion label) \label{MS_sigma} If $\theta\in 2\mathbb Z+1$, then every projective equivalence class of spherical tori, with a labelled $2$-torsion point, can be parametrized by $\mathbb R$. In each class there is a unique surface having the isometric spherical triangle decomposition as above, with the labelled $2$-torsion lying on $L_1$. Moreover, the three angles of the two spherical triangles must be integral multiples of $\pi$.\\

\begin{figure}
\includegraphics[scale=0.7]{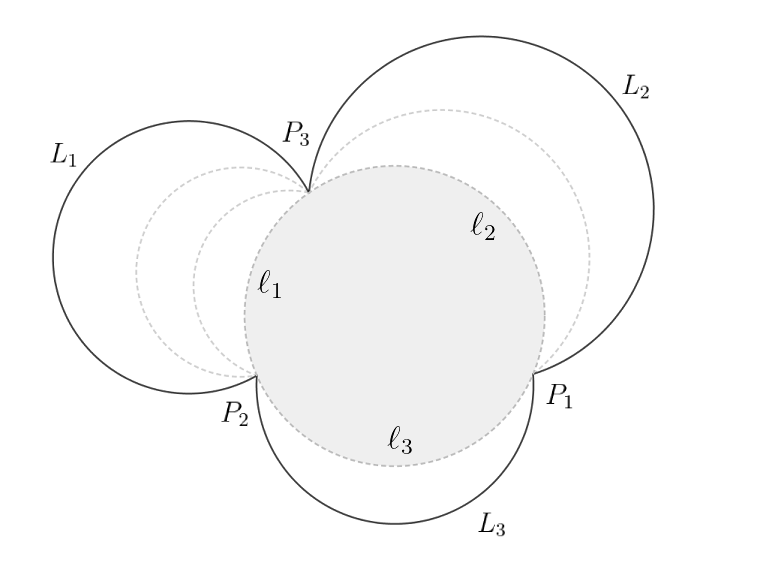}\\
\caption{A demonstration of a spherical triangle with interior angles $(\pi\theta_1,\pi\theta_2,\pi\theta_3)=(4\pi,5\pi,6\pi)$. All the bounded regions in the figure are hemispheres, and the shaded region is the basic spherical triangle.}
\end{figure}

We give a detailed description for such spherical tori. Let $\pi\theta_1,\pi\theta_2,\pi\theta_3$ be the three interior angles of the two spherical triangles, so $\theta=\theta_1+\theta_2+\theta_3=2n+1$, with $\theta_i$ all integers and satisfying the triangle inequalities. The spherical triangle can be obtained from contiguously gluing hemispheres along three edges of a basic spherical triangle of interior angle $\pi,\pi,\pi$ (thus also a hemisphere), see Figure 1.\\  

As a consequence, the set of spherical tori of angle $2n+1$, with a labelled 2-torsion, is parametrized by \begin{align*}\{(\theta_1,\theta_2,\theta_3)\in\mathbb \{1,\dots,n\}^3, \theta_1+\theta_2+\theta_3=2n+1\}\\\times\{(\ell_1,\ell_2,\ell_3)\in (\mathbb R^+)^3, \ell_1+\ell_2+\ell_3=2\pi\}\times \mathbb R,\end{align*} where $\ell_1,\ell_2,\ell_3$ are parameters for the lengths of edges of the basic spherical triangle. Quotiening the last component and the $\mathbb Z_3$-action on cyclically permuting $\theta_i$ and $\ell_i$ gives the set of Lam\'e equations with unitary monodromy for $n\in\mathbb Z_{>0}$.

\section{Main Theorems}

In this section we will prove the two conjectures of Dahmen. We first establish the relation between the monodromy group and the shape of the spherical triangle.

\definition Denote by 
$S_{\theta_1,\theta_2,\theta_3}(\ell_1,\ell_2,\ell_3)$ the spherical torus with decomposition as in Theorem \ref{MS_sigma} and with those parameters, and with a labelled 2-torsion at the midpoint of $L_1$.\\

It is noteworthy that $f$ maps the boundaries of every hemispheres in $\triangle$ and $\triangle'$ to the unit circle of $\mathbb{CP}^1$, and the centers of these hemispheres are mapped to $0$ or $\infty$, regardless of the monodromy. As a result, these centers are $a_i$ or $-a_i$ in the ansatz solution respectively. One simple observation is

\proposition \label{dev_mono} The monodromy of $ S_{\theta_1,\theta_2,\theta_3}(\ell_1,\ell_2,\ell_3)$ is given by $f(z+\omega_1)=e^{\mp i(\ell_2+\ell_3)}f(z)$ and $f(z+\omega_2)=e^{\pm i(\ell_1+\ell_3)}f(z)$. 
\proof Since attaching an even number of hemispheres on the sides does not affect the projective monodromy, we may assume $\triangle$ and $\triangle'$ are basic (that is, $n=1$). We first consider the case the center of $\triangle$ is mapped to $0$. We glue $\triangle$ and $\triangle'$ along $L_3$, and we find the monodromy of the loop starting from any point on $L_2$ to the corresponding point on $L_2'$. Geometrically we can see that the loop is homotopic to a path on unit circle, clockwise of length $\ell_2+\ell_3$, thus $f(z+\omega_1)=e^{-i(\ell_2+\ell_3)}f(z)$ for $z$ on the unit circle. Similarly $f(z+\omega_2)=e^{i(\ell_1+\ell_3)}f(z)$, and they hold for any $z\in \mathbb{CP}^1$ since the monodromy of $f$ is in $\mathrm{PSU}(2)$. If the center of $\triangle$ is mapped to $\infty$, then $f$ has inverse monodromy. \qedhere\\

\begin{figure}\label{sphericalTrianglePic}
\includegraphics[scale=0.4]{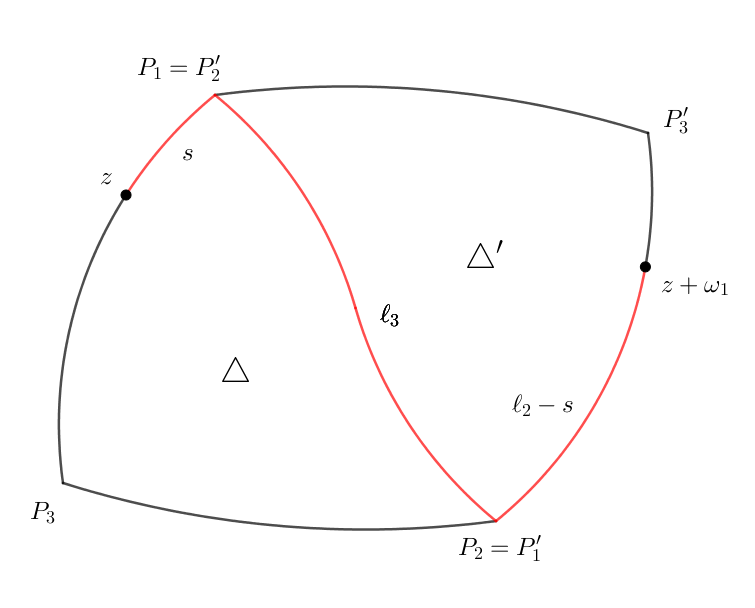}\\
\caption{The $n=1$ case. The image of the red path under the map $f$ is an arc of length $\ell_2+\ell_3$ on the unit circle.}
\end{figure}

From this we can easily deduce the constraint on projective monodromy.
\proposition \label{proj_mono} The projective monodromy parameters $2s,2t(\mathrm{mod}\ 1)$ must satisfies the restriction $2s\neq 0,2t\neq 0, 2s+2t\neq 0(\mathrm{mod}\ 1)$. Any such $2s,2t$ would give rise to $n(n+1)/2$ Lam\'e equations with unitary monodromy satisfying the monodromy (\ref{proj_mono_eq}). The distribution of the monodromy parameters is shown in Figure 3.
\proof We show that if the center of $\triangle$ is mapped to $0$, then we have \[2s<1,2t<1\textnormal{ and }2s+2t>1.\] Comparing Proposition \ref{dev_mono} and (\ref{proj_mono_eq}), we have \[2s=(\ell_2+\ell_3)/2\pi\textnormal{ and }2t=(\ell_1+\ell_3)/2\pi.\]The inequalities then follows from \[\ell_2+\ell_3<2\pi, \ell_1+\ell_3<2\pi\textnormal{ and }\ell_1+2\ell_2+\ell_3=2\pi+\ell_2>2\pi.\] Conversely, if $2s$ and $2t$ satisfies the inequalities, then we can solve \[(\ell_1,\ell_2,\ell_3)=(2\pi-4\pi t,2\pi-4\pi s,4\pi t+4\pi s-2\pi),\]and we can find one projective equivalence class in each connected component of $MS_{1,1}^{[2]}(2n+1)$ with such monodromy. If the center of $\triangle$ is mapped to $\infty$, then $f$ has inverse monodromy, and we have  \[2s>0,2t>0\textnormal{ and }2s+2t<1.\] \qedhere\\

For the ordinary monodromy, we need to take the attached hemispheres on the sides into consideration. 
\proposition \label{ordi_mono} For fixed parameters $\theta_1,\theta_2,\theta_3\in\{1,2,\dots,n\}$ with $\theta=\theta_1+\theta_2+\theta_3=2n+1$, the monodromy parameter $s,t (\mathrm{mod}\ 1)$ for the ansatz $w_a$ satisfies either \[s<\frac{\theta_1}{2},\ t<\frac{\theta_2}{2},\ s+t>\frac{\theta_1+\theta_2-1}{2}\] or \[s>-\frac{\theta_1}{2},\ t>-\frac{\theta_2}{2},\ s+t<-\frac{\theta_1+\theta_2-1}{2}.\] 

Translating these regions into $(0,1)\times (0,1)$, we have the distribution of monodromy parameters shown in Figure 4.
\proof Note that for any fixed $(2s,2t)$, there are $4$ choices of $(s,t)$, so we ought to determine for each $\theta_1,\theta_2,\theta_3$ which of the $8$ regions in Figure 4 $(s,t)$ should lies in. For the case $n=1$ and $\theta_1=\theta_2=\theta_3=1$, it is known that $(s,t)$ lies in the region we have described (\cite{A},p30). For $n>1$, note that attaching two hemispheres on the side of $L_2$ or $L_3$ will translate the parameter $s \mathrm{(mod}\ 1)$ by $\frac{1}{2}$, so in total $s$ is translated by $\frac{\theta_1-1}{2}$. A similar argument holds for $t$. \qedhere\\

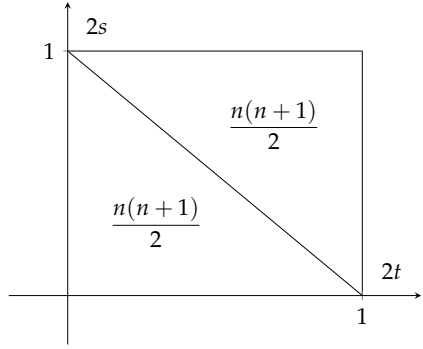
\begin{figure}
\begin{tikzpicture}[scale=0.8]
\begin{axis}[
    xmin=-0.2, xmax=1.2,
    xtick={0,1},
    ymin=-0.2, ymax=1.2,
    ytick={0,1},
    axis lines=center,
    axis on top=true,
    domain=0:1,
    ]
    \addplot [mark=none] coordinates {(1, 0) (0, 1)};
    \addplot [mark=none] coordinates {(0, 1) (1, 1)};
    \addplot [mark=none] coordinates {(1, 0) (1, 1)};
    \node at (axis cs: 0.7,0.7) {$\dfrac{n(n+1)}{2}$};
    \node at (axis cs: 0.3,0.3) {$\dfrac{n(n+1)}{2}$};
    \node at (axis cs: 1.1,0.1) {$2t$};
    \node at (axis cs: 0.1,1.1) {$2s$};
\end{axis}
\end{tikzpicture}
\caption{The number of developing maps with unitary monodromy and with given parameter $2s,2t (\mathrm{mod}\ 1)$.}
\end{figure}

\begin{figure}
\begin{subfigure}{0.45\textwidth}
\begin{tikzpicture}[scale=0.8]
\begin{axis}[
    xmin=-0.2, xmax=1.2,
    xtick={0,0.5,1},
    ymin=-0.2, ymax=1.2,
    ytick={0,0.5,1},
    axis lines=center,
    axis on top=true,
    domain=0:1,
    ]
    \addplot [mark=none] coordinates {(0.5, 0) (0, 0.5)};
    \addplot [mark=none] coordinates {(1, 0) (0, 1)};
    \addplot [mark=none] coordinates {(1, 0.5) (0.5, 1)};
    \addplot [mark=none] coordinates {(0, 0.5) (1, 0.5)};
    \addplot [mark=none] coordinates {(0, 1) (1, 1)};
    \addplot [mark=none] coordinates {(0.5, 0) (0.5, 1)};
    \addplot [mark=none] coordinates {(1, 0) (1, 1)};
    \node at (axis cs: 0.15,0.12) {$\dfrac{l(l-1)}{2}$};
    \node at (axis cs: 0.65,0.12) {$\dfrac{l(l-1)}{2}$};
    \node at (axis cs: 0.15,0.62) {$\dfrac{l(l-1)}{2}$};
    \node at (axis cs: 0.65,0.62) {$\dfrac{l(l+1)}{2}$};
    \node at (axis cs: 0.35,0.37) {$\dfrac{l(l+1)}{2}$};
    \node at (axis cs: 0.85,0.37) {$\dfrac{l(l-1)}{2}$};
    \node at (axis cs: 0.35,0.87) {$\dfrac{l(l-1)}{2}$};
    \node at (axis cs: 0.85,0.87) {$\dfrac{l(l-1)}{2}$};
    \node at (axis cs: 1.1,0.1) {$t$};
    \node at (axis cs: 0.1,1.1) {$s$};
\end{axis}
\end{tikzpicture}
\end{subfigure}
\begin{subfigure}{0.45\textwidth}
\begin{tikzpicture}[scale=0.8]
\begin{axis}[
    xmin=-0.2, xmax=1.2,
    xtick={0,0.5,1},
    ymin=-0.2, ymax=1.2,
    ytick={0,0.5,1},
    axis lines=center,
    axis on top=true,
    domain=0:1,
    ]
    \addplot [mark=none] coordinates {(0.5, 0) (0, 0.5)};
    \addplot [mark=none] coordinates {(1, 0) (0, 1)};
    \addplot [mark=none] coordinates {(1, 0.5) (0.5, 1)};
    \addplot [mark=none] coordinates {(0, 0.5) (1, 0.5)};
    \addplot [mark=none] coordinates {(0, 1) (1, 1)};
    \addplot [mark=none] coordinates {(0.5, 0) (0.5, 1)};
    \addplot [mark=none] coordinates {(1, 0) (1, 1)};
    \node at (axis cs: 0.15,0.12) {$\dfrac{l(l+1)}{2}$};
    \node at (axis cs: 0.65,0.12) {$\dfrac{l(l+1)}{2}$};
    \node at (axis cs: 0.15,0.62) {$\dfrac{l(l+1)}{2}$};
    \node at (axis cs: 0.65,0.62) {$\dfrac{l(l-1)}{2}$};
    \node at (axis cs: 0.35,0.37) {$\dfrac{l(l-1)}{2}$};
    \node at (axis cs: 0.85,0.37) {$\dfrac{l(l+1)}{2}$};
    \node at (axis cs: 0.35,0.87) {$\dfrac{l(l+1)}{2}$};
    \node at (axis cs: 0.85,0.87) {$\dfrac{l(l+1)}{2}$};
    \node at (axis cs: 1.1,0.1) {$t$};
    \node at (axis cs: 0.1,1.1) {$s$};
\end{axis}
\end{tikzpicture}
\end{subfigure}
\caption{The number of ansatz solutions with unitary monodromy and with given parameter $s,t (\mathrm{mod}\ 1)$ and $n=2l-1$ (left) or $n=2l$ (right).}
\end{figure}
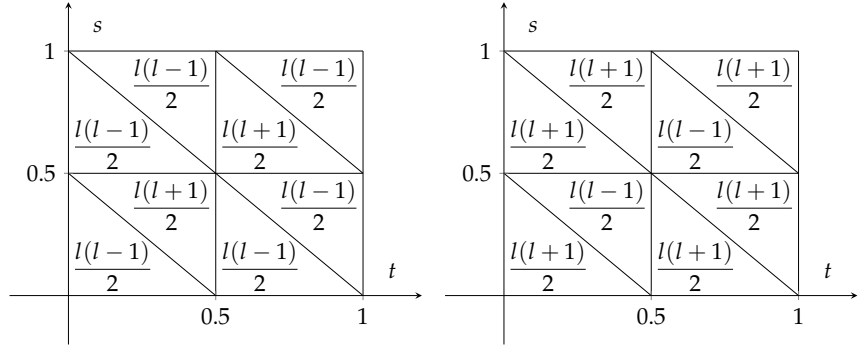

Now we are ready to prove Dahmen's conjectures with the propositions above.\\

\paragraph{\it{Proof of Main Conjectures}} Note that if $n=3$ and $3|(N-1)$, then there is a unique spherical torus $S_{\frac{N-1}{3},\frac{N-1}{3},\frac{N-1}{3}}(2\pi/3,2\pi/3,2\pi/3)$ fixed by the $\mathbb Z/3$-action on the labels, thus $3L_n(N)-2\epsilon(n,N)$ counts the number of Lam\'e equations with parameter \[(2s,2t)=(\dfrac{k_1}{N},\dfrac{k_2}{N}), 0<k_1<N-k_2<N, \textrm{gcd}(k_1,k_2,N)=1\textnormal{ and }N\geq 3.\] 

For convenience we suppose this also holds for $N=1,2$. By counting lattice points and Proposition \ref{proj_mono} we have
\[\sum_{d|N}(3L_n(d)-2\epsilon(n,d))=\dfrac{n(n+1)}{2}(N^2-3N+2).\]
The formula (\ref{fin_proj_mono}) then follows from M\"obius inversion.\\

Similarly we have that $3L_n'(N)-2\epsilon(n,N)$ counts the number of Lam\'e equations with parameter \[(s,t)=(\dfrac{k_1}{N},\dfrac{k_2}{N}), 0<k_1<N-k_2<N,  \mathrm{gcd}(k_1,k_2,N)=1\textnormal{ and }N\geq 3,\] and we suppose this holds for $N=1,2$. By counting lattice points and Proposition \ref{ordi_mono} we have
\begin{align*}&\sum_{d|N}(3L_n'(d)-2\epsilon(n,d))\\&=\begin{cases}
a_n\dfrac{3(m-1)(m-2)}{2}+(b_n-a_n)\dfrac{m(m-1)}{2} \qquad\textnormal{if }N=2m-1,\\
a_n\dfrac{3(m-1)(m-2)}{2}+(b_n-a_n)\dfrac{(m-1)(m-2)}{2} \qquad\textnormal{if }N=2m,\end{cases}
\\&=\begin{cases}
\dfrac{n(n+1)}{2}\dfrac{m(m-1)}{2}-3a_nm\qquad\textnormal{if }N=2m-1,\\
\dfrac{n(n+1)}{2}\dfrac{m^2}{2}-(2a_n+b_n)\dfrac{3m-2}{2}\qquad\textnormal{if }N=2m,
\end{cases}
\\&=\begin{cases}
\dfrac{n(n+1)}{16}(N^2-1)-\dfrac{3}{2} a_n (N+1) \qquad\textnormal{if }N=2m-1,\\
\dfrac{n(n+1)}{16}N^2-(2a_n+b_n)(\dfrac{3N}{4}-1) \qquad\textnormal{if }N=2m.\end{cases}
\end{align*}

The formula (\ref{fin_ordi_mono}) then follows from M\"obius inversion. The proof is complete. \qedhere\\

\begin{paragraph}{\textbf{Remark}} We can compare our proof of Conjecture 1.1 with Dahmen's proof using dessin d'enfant. 

Given a spherical torus  $(S,p)=S_{\theta_1,\theta_2,\theta_3}(\ell_1,\ell_2,\ell_3)$ with $\ell_i=2\pi\dfrac{m_i}{N}, i=1,2,3$ and $m_i$ positive integers, we assume the developing map $f$ sends $p$ to some $N$th root of unity, and let \[g(z)=\left(\dfrac{1-z^N}{1+z^N}\right)^2, h(z)=\dfrac{z}{z-1}.\]

The composition $h\circ g\circ f$ sends the $2N$-division points on the boundaries of hemispheres to $0$ or $1$ alternately, and also sends the centers of the hemispheres to $\infty$. As $h\circ g \circ f$ is independent of monodromy and ramify only at $0,1,\infty$, it serves as a Belyi function for the underlying elliptic curve of $S$. The dessin d'enfant corresponding to this Belyi function consists of $2n+1$ loops of $3$ directions through $p$ on the elliptic curve. The $i$th direction has $2n+1-2\theta_i$ loops with the numbers of edges alternate between $2m_i$ and $2N-2m_i$. Taking the quotient of the elliptic curve as well as the dessin by the involution $z\mapsto -z$, we obtain a dessin on $\mathbb P^1$ of Type I introduced in Dahmen's proof of Conjecture 1.1. (See Figure 5.). As a consequence, counting Type I dessins on $\mathbb P^1$ is equivalent to counting spherical tori of finite monodromy we have done here.

\begin{figure}

\begin{subfigure}{0.6\textwidth}
\includegraphics[scale=0.3]{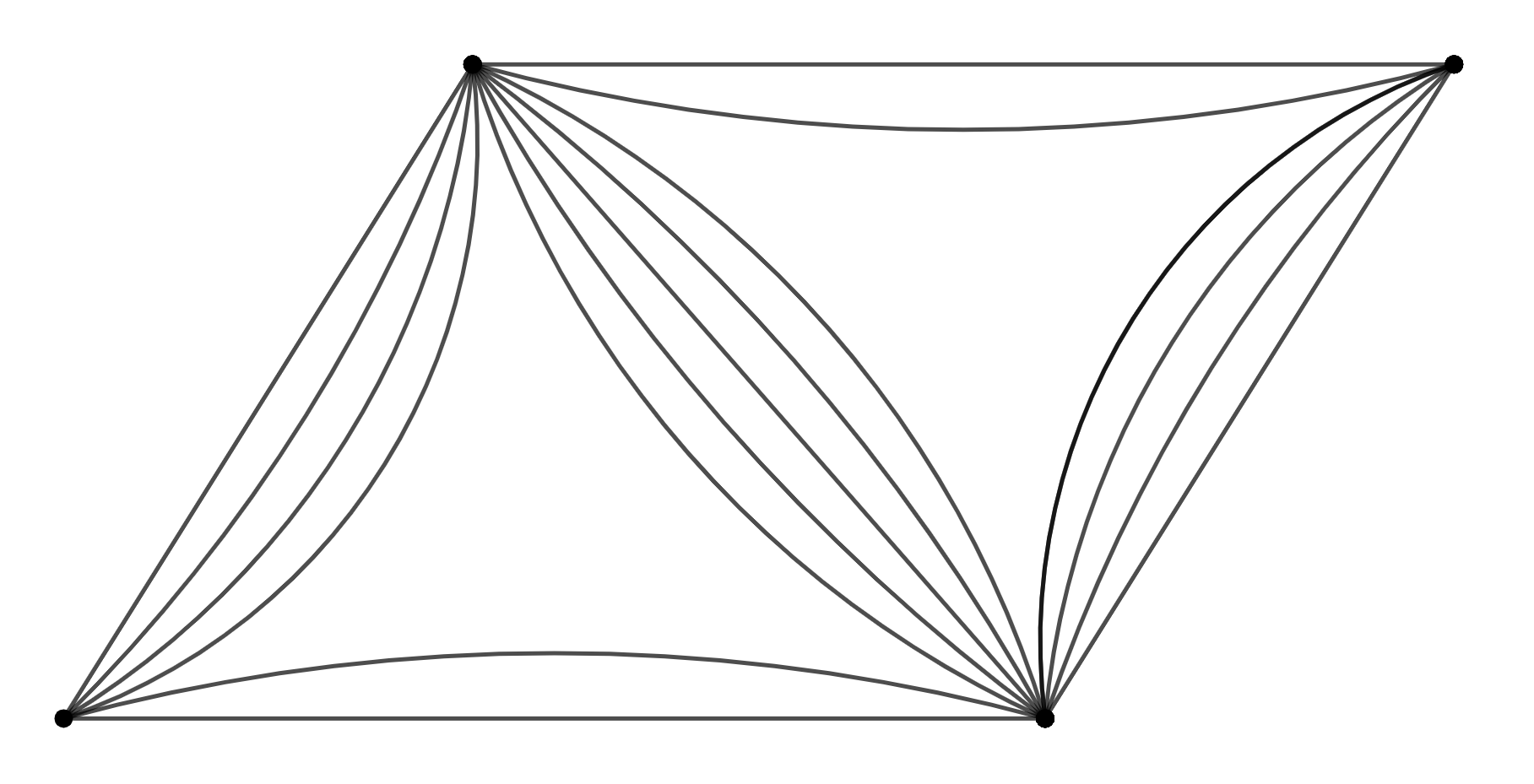}\\
\end{subfigure}
\begin{subfigure}{0.3\textwidth}
\includegraphics[scale=0.3]{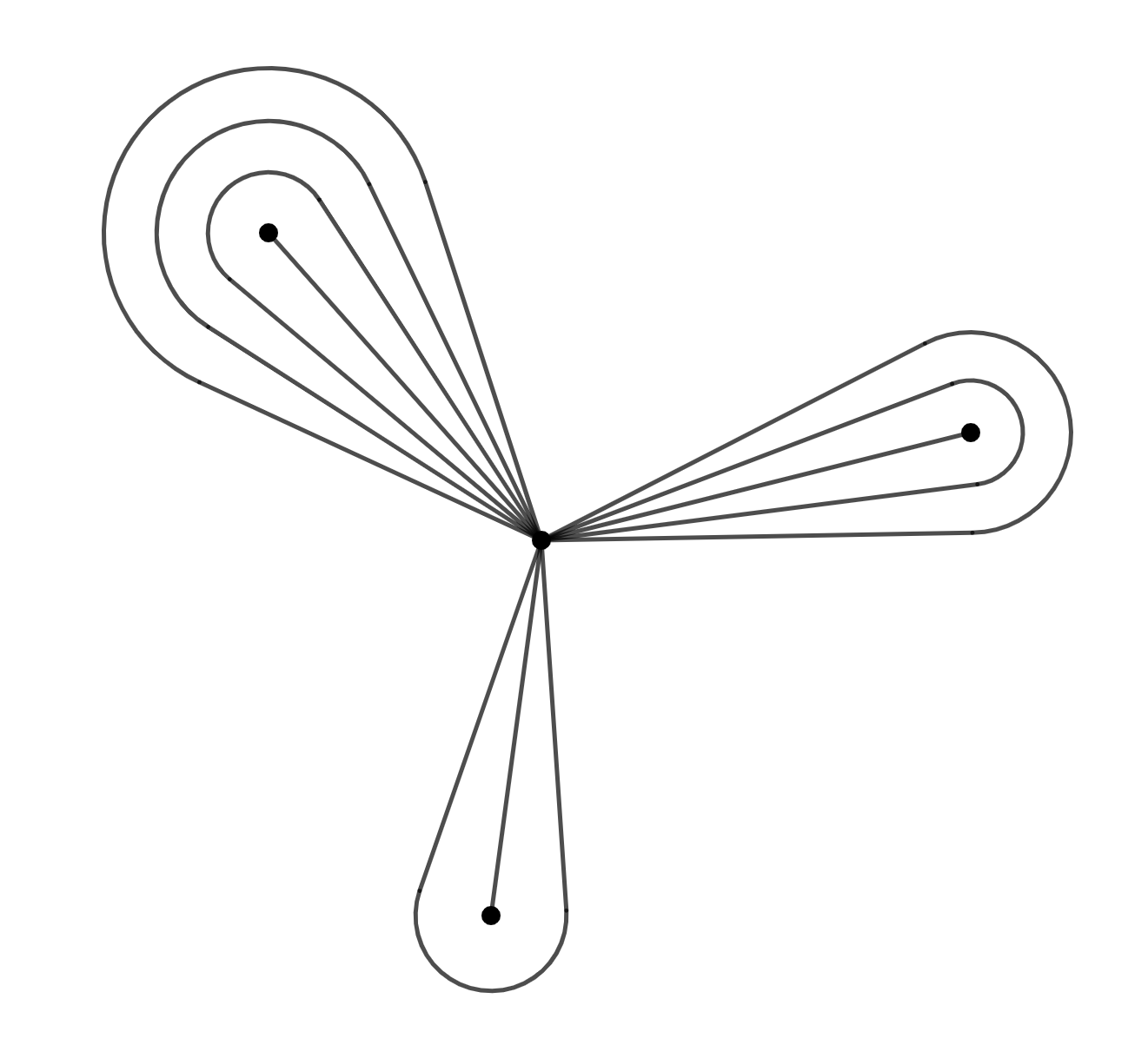}\\
\end{subfigure}

\caption{The dessin corresponding to the spherical torus obtained from Figure 2 (the opposite sides of the parallelogram are identified) (left), and its quotient dessin on $\mathbb P^1$ (right). Note that each segment represents several edges of the dessin (according to $m_i$).}
\end{figure}

\end{paragraph}

\end{document}